\documentclass[12pt,leqno]{article}
\usepackage{amsfonts}
\usepackage{ccaption,mathrsfs,amsmath,amssymb}
\usepackage{amsmath, amscd, amsfonts, amsthm, tikz, times}

\newcounter{conjecture}
\newcounter{remark}
\newcounter{corollary}
\newenvironment{conjecture}{\medskip{\bf Conjecture.\ }\em}{\rm}

\newcommand{\eqnsection}{
    \renewcommand{\theequation}{\thesection.\arabic{equation}}
    \makeatletter
    \csname @addtoreset\endcsname{equation}{section}
    \makeatother}
\newtheorem{theorem}{Theorem}
\newtheorem{lemma}{Lemma}

\newtheorem{prop}{Proposition}

\newcommand{\dd}{\delta}
\newcommand{\DD}{\Delta}

\newcommand{\lar}{\longrightarrow}

\newcommand{\om}{\Omega}
\newcommand{\aaa}{\alpha}

\newcommand{\CC}{\mathbb{C}}

\def \be{\begin{equation}}
\def \ee{\end{equation}}
\def \bt{\begin{theorem}}
\def \et{\end{theorem}}
\def \bea{\begin{eqnarray}}
\def \eea{\end{eqnarray}}
\def \bas{\begin{eqnarray*}}
\def \eas{\end{eqnarray*}}

\newcommand {\rrr}[1]{(\ref{#1})}

\def \Be{\beta}
\def \ga{\gamma}
\def \Ga{\Gamma}

\def \om{\omega}

\def \si{\sigma}

\def \th{\theta}

\def \ze{\zeta}

\def \ff{\infty}

\def \DD{\mathbb{D}}

\def \HH{{\mathbb H}}

\def \KK{{\cal K}}

\def \RR{{\mathbb R}}
\def \SS{{\cal S}}

\def \({\left(}
\def \){\right)}

\def \vski{\vspace{12pt}}

\def \bc{\begin{center} }
\def \ec{\end{center} }
\def \bs{\begin{slide} }
\def \es{\end{slide} }

\def\square{{\vcenter{\vbox{\hrule height.3pt
         \hbox{\vrule width.3pt height5pt \kern5pt
            \vrule width.3pt}
         \hrule height.3pt}}}}
\def\qed{{\hfill $\Box$ \bigskip}}

\newcounter{cccases}
\newcommand{\ccases}[1]{\begingroup \refstepcounter{cccases} {\bf \fontsize{12}{12}\selectfont Example \thecccases }  \label{#1}\endgroup}

\eqnsection
\begin{document}

\title{On the expected exit time of planar Brownian motion from simply connected domains.}

\author{
\begin{tabular}{c}
\textit{Greg Markowsky} \\
gmarkowsky@gmail.com \\
+61 03 9905 4487 \\
Monash University \\
%Department of Mathematical Sciences \\
Victoria, 3800 Australia
\end{tabular}}

\bibliographystyle{amsplain}
\maketitle \eqnsection \setlength{\unitlength}{2mm}

\begin{abstract}
This paper presents some results on the expected exit time of Brownian motion from simply connected domains in $\CC$. We indicate a way in which Brownian motion sees the identity function and the Koebe function as the smallest and largest analytic functions, respectively, in the Schlicht class. We also give a sharpening of a result of McConnell's concerning the moments of exit times of Schlicht domains. We then show how a simple formula for expected exit time can be applied in a series of examples. Included in the examples given are the expected exit times from given points of a cardioid and regular $m$-gon, as well as bounds on the expected exit time of an infinite wedge. We also calculate the expected exit time of an infinite strip, and in the process obtain a probabilistic derivation of Euler's result that $\zeta(2)=\sum_{n=1}^\ff \frac{1}{n^2}= \frac{\pi^2}{6}$. We conclude by showing how the formula can be applied to some domains which are not simply connected.
\end{abstract}

\section{Introduction and statement of results}

A domain $U$ in $\CC$ is {\it simply connected} if any curve in $U$ can be homotopically deformed to a point. Equivalently, $U$ is simply connected if $U^c$ is connected, where the complement is taken in the Riemann sphere, $\hat{\CC}=\CC \cup \{\ff\}$, which is the one point compactification of $\CC$. Such domains are of paramount importance in complex analysis, due in large part to the Riemann Mapping Theorem, which states that for any simply connected domain $U \subsetneq \CC$ and point $z \in U$ there is a conformal map from $\DD$ onto $U$ with $f(0)=z$. The {\it Schlicht class} $\SS$ is the set of all conformal functions $f$ on $\DD$ normalized so that $f(0)=0, f'(0)=1$. If a domain $U$ is equal to $f(\DD)$ for some $f\in \SS$, we will call $U$ a {\it Schlicht domain}. Clearly, if $U$ is any simply connected domain smaller than $\CC$ itself and $z \in U$ we may find a linear map sending $z$ to 0 and $U$ to a Schlicht domain. In this way the study of $\SS$ is no less general than the study of arbitrary simply connected domains. In a number of different ways, the identity function $I(z)=z$ with image $I(\DD)=\DD$ is considered to be the smallest in $\SS$, and the Koebe function $K(z) = \frac{z}{(1-z)^2}$ with image $K(\DD)=\CC \backslash (\ff, -1/4]$ is considered to be the largest. We may view this intuitively as follows. The normalization creates an equivalence between the notion of "size" and the notion of "closeness of $0$ to the boundary". The disc is the only domain in which $0$ is equally close to every boundary point. In that sense, $0$ is closer to the boundary in the unit disc than in any other Schlicht domain. On the other hand, the smallest nonempty complement of a simply connected domain seems likely to be a half-line $(-\ff,0]$, and a point $x$ on $(0,\ff)$ should be farther from the boundary than $e^{i \th} x$, $\th \in (0,\pi)$, for instance. Arguing in this manner indicates that $K(\DD) = \CC \backslash (\ff, -1/4]$ should be the Schlicht domain for which $0$ is farthest from the boundary. This intuition may be useful as far as it goes, but a far more sophisticated statement is the following celebrated theorem of de Branges, which was originally conjectured by Bieberbach.

\begin{theorem} \label{deb}
If $f \in \SS$ with $f(z) = z + \sum_{n=2}^{\ff} a_n z^n$, then $|a_n| \leq n$ for all $n \geq 2$. If for any such $n$ we have $|a_n|=n$, then $f$ is of the form $e^{-i \aaa} K(e^{i \aaa}z)$ for real $\aaa$.
\end{theorem}

Thus, out of all functions in $\SS$, the Koebe function and its rotations have the largest derivatives at $0$. The highly nontrivial proof of Theorem \ref{deb}, which follows decades of work by analysts, can be found in \cite{deb}.

\vski

These concepts can be connected to probabilistic considerations in the following way. The expected amount of time that it takes a Brownian motion starting at a point $z$ to leave a domain $U$ gives some sort of measure of the size of $U$ and the distance between $z$ and $U^c$. This was perhaps first noted by Davis in \cite{davis}. In that paper the question was raised of finding a way in which Brownian motion identifies the Koebe domain as the largest Schlicht domain and the unit disk as the smallest. Initially we may hope that the expected amount of time that it takes for a Brownian motion starting at $0$ to leave a Schlicht domain is bounded above and below by the corresponding times for the Koebe domain and unit disc, respectively. This turns out to be correct for the unit disc, and we may say more. In \cite{mac}, McConnell proved the following.

\begin{theorem} \label{}
Let $\Phi$ be a nonnegative, nondecreasing convex or concave function, where if $f$ is concave then we also require that $\Phi(e^{2x})$ is convex. Then, for any $f \in \SS$ we have

\begin{equation} \label{}
E_0[\Phi(\tau(\DD))] \leq E_0[\Phi(\tau(f(\DD)))]
\end{equation}

In particular, $E_0[\tau(\DD)^p] \leq E_0[\tau(f(\DD))^p]$ for $0<p<\ff$.
\end{theorem}

We will prove the following addition to this result.

\begin{theorem} \label{add}
Let $\Phi$ be a nonnegative, strictly increasing convex or strictly concave function, where if $f$ is strictly concave then we also require that $\Phi(e^{2x})$ is convex. Assume further that $E_0[\Phi(\tau(\DD))]<\ff$. Then, for any $f \in \SS, f \neq I$ we have

\begin{equation} \label{}
E_0[\Phi(\tau(\DD))] < E_0[\Phi(\tau(f(\DD)))]
\end{equation}

In particular, $E_0[\tau(\DD)^p] < E_0[\tau(f(\DD))^p]$ for $0<p<\ff$.
\end{theorem}

Such a statement for $K(\DD)$ is harder to formulate and prove. To begin with, many of the moments of $\tau(K(\DD),0)$ are infinite. It is a consequence of \cite[Thm. 4.1]{burk} that $E[\tau(K(\DD),0)^p]<\ff$ if and only if $p<1/4$, and it seems likely that $E[\tau(K(\DD),0)^p] > E[\tau(f(\DD),0)^p]$ for all $p<1/4$ and $f \in \SS \backslash \KK$, where we let $\KK$ denote the Koebe function and its rotations, that is, $\KK = \{e^{-i \aaa} K(e^{i \aaa}z) : \aaa \in {\mathbb R}\}$. Unfortunately, this statement, if true, remains unproved. We consider a different approach to the problem for $K(\DD)$. We will use $f(r\DD)$ as a bounded approximation to $f(\DD)$ for all $f \in \SS$. This is justified by the fact that $f(r\DD) \lar f(\DD)$ in the sense of Carath\'eodory as $r \lar 1$(see \cite[Ch. 3]{dur}), and the approximation is in some sense uniform between domains. We will then show that Brownian motion leaves $K(r\DD)$ more quickly than $f(r\DD)$, where $f$ is any function in $\SS \backslash \KK$.

\begin{theorem} \label{greg}
If $f \in \SS$ and $r \leq 1$, then

\begin{equation} \label{}
E_{0}[\tau(I(r\DD))] \leq E_{0}[\tau(f(r\DD))] \leq E_{0}[\tau(K(r\DD))]
\end{equation}

If $E_{0}[\tau(f(r\DD))] = E_{0}[\tau(I(r\DD))]$ for any $r \leq 1$ then $f(z)=I(z)=z$, and if $E_{0}[\tau(f(r\DD))] = E_{0}[\tau(K(r\DD))]$ for any $r<1$ then $f\in \KK$.
\end{theorem}

In this sense, Brownian motion sees $K(r\DD)$ as larger than $f(r\DD)$ for any $f \in \SS \backslash \KK$. We will also show

\begin{theorem} \label{kris}
If $f \in \SS, f \notin \KK$, then $E_{0}[\tau(K(r\DD))] - E_{0}[\tau(f(r\DD))]$ is an increasing function of $r$ which is positive for $r>0$.
\end{theorem}

We may intuitively rephrase this result as "Brownian motion would take longer to exit $K(\DD)$ than any other Schlicht domain were it not that many of the expected exit times are infinite".

\section{Proofs of Theorems \ref{add}, \ref{greg}, and \ref{kris}}

Any analytic function can be expressed locally as a Taylor series(see \cite{rud}, \cite{ahl}, or \cite{schaum}). The Parseval Identity for analytic functions(see \cite[Thm. 10.22]{rud}) states that

\begin{equation} \label{par}
\frac{1}{2\pi} \int_{0}^{2\pi}|g(s e^{i \th})|^2 d\th = \sum_{n=0}^{\ff} |b_n|^2 s^{2n} ,
\end{equation}

for any $g(z) = \sum_{n=0}^{\ff} b_n z^n$ analytic on an open set containing $s\DD$. This equality will be the key for the proof of all three theorems.

\subsection*{Proof of Theorem \ref{add}}

Assume first that $\Phi$ is convex. We will follow the method given in \cite{mac}. Assume $f \in \SS, f \neq I$. $f(B_t)$ is a time-changed Brownian motion(see \cite{durBM} or \cite{revyor}). That is, there is a complex Brownian motion $\hat{B}_t$ such that $f(B_t) = \hat{B}_{\si_t}$, where $\si_t = \int_{0}^{t}|f'(B_s)|^2 ds$. It follows that

\begin{equation} \label{}
E_{0}[\tau(f(\DD))] = E \int_{0}^{\tau(\DD,0)} |f'(B_s)|^2 ds
\end{equation}

The process $e^{i\th}B_s$ is again a planar Brownian motion for real $\th$. Let $f(z)$ have the Taylor expansion $\sum_{n=1}^\ff a_n z^n$. The assumptions on $f$ imply $a_1=1, a_m \neq 0$ for some $m\geq 2$. We therefore have by Jensen's inequality

\begin{equation} \label{}
\begin{split}
E_0[\Phi(\tau(f(\DD)))] &= \frac{1}{2\pi} \int_{0}^{2\pi} E \Phi \Big( \int_0^{\tau(\DD,0)} |f'(e^{i\th}B_s)|^2 ds\Big)d\th \\
& \geq E \Phi \Big( \int_0^{\tau(\DD,0)}\frac{1}{2\pi} \int_{0}^{2\pi} |f'(e^{i\th}B_s)|^2 d\th ds \Big) \\
& = E \Phi \Big( \int_0^{\tau(\DD,0)} \sum_{n=1}^{\ff} n^2 |a_n|^2 |B_s|^{2n-2} ds\Big) \\
& > E \Phi \Big( \int_0^{\tau(\DD,0)} |a_1|^2 ds \Big) \\
& = E \Phi \Big( \int_0^{\tau(\DD,0)} \! \! \! \! \! ds \Big) = E \Phi(\tau(\DD,0)) = E_0[\Phi(\tau(\DD))]
\end{split}
\end{equation}

where Parseval's Identity \rrr{par} was applied to $f'(z) = \sum_{n=1}^\ff n a_n z^{n-1}$.

\vski

Suppose now that $\Phi$ is strictly concave and $\Phi(e^{2x})$ is convex. Examining the proof corresponding to this statement in \cite{mac}, we see that the following consequence of Jensen's inequality for concave functions was used.

\begin{equation} \label{}
\Phi \Big( \int_{0}^{\tau} |f'(e^{i \th}B_s)|^2 ds \Big) \geq \int_{0}^{\tau} \Phi(\tau |f'(e^{i \th}B_s)|^2) \frac{ds}{\tau}
\end{equation}

However, with strictly concave $\Phi$, Jensen's inequality is easily seen to be strict when $\tau|f'(e^{i \th}B_s)|^2$ is nonconstant and the integrals in question are finite. We see that if $E_0[\Phi(\tau(f(\DD)))]< \ff$ then McConnell's proof in \cite{mac} is sufficient to prove our result. Clearly, if $E_0[\Phi(\tau(f(\DD)))]= \ff$ then the result holds as well.

\vski

It remains only to show that $E[\tau(\DD)^p] < \ff$ for all $p>0$. There must surely be many proofs of this fact, but we will be content to state it as an immediate consequence of \cite[Thm. 4.1]{burk}. \qed

\subsection*{Proof of Theorems \ref{greg} and \ref{kris}}

Both theorems are a direct consequence of Lemma \ref{bigguy} below. Due to the simplicity of the proof, it seems almost certain that this lemma has been noted before. Nevertheless, a literature search has failed to yield a suitable reference. We therefore include a proof.

\begin{lemma} \label{bigguy}
Suppose $f(z) = \sum_{n=0}^{\ff} a_n z^n$ is conformal on $\DD$. Then, for $0<r \leq 1$ we have

\begin{equation} \label{re}
E_{f(0)}[\tau(f(r\DD))] = \frac{1}{2}\sum_{n=1}^{\ff} |a_n|^2 r^{2n}
\end{equation}

In particular,

\begin{equation} \label{re2}
E_{f(0)}[\tau(f(\DD))] = \frac{1}{2}\sum_{n=1}^{\ff} |a_n|^2
\end{equation}

\end{lemma}

{\bf Proof of lemma:} Assume first that $r<1$. As before, due to the fact that $f(B_t)$ is a time-changed Brownian motion we have

\begin{equation} \label{}
E_{f(0)}[\tau(f(r\DD))] = E \int_{0}^{\tau(r\DD,0)} |f'(B_s)|^2 ds
\end{equation}

Set $\si_{t} = \int_{0}^{t} |f'(B_s)|^2 ds$. Applying the Optional Stopping Theorem to the martingale $|f(B_t)|^2 - 2 \si_t$ gives

\begin{equation} \label{}
\begin{split}
2 E_{f(0)}[\tau(f(r\DD))] & = 2 E[\si_{\tau(r\DD,0)}] \\
& = E[|f(B_{\tau(r\DD,0)})|^2] \\ &= \frac{1}{2\pi} \int_{0}^{2\pi} |f(re^{i\th})|^2 d\th \\ &= \sum_{n=1}^{\ff} |a_n|^2r^{2n}
\end{split}
\end{equation}

In order to set $r=1$ to obtain the final statement, we note first that if $\sum_{n=1}^{\ff} |a_n|^2 = \ff$ then $E_{f(0)}[\tau(f(r\DD))] \lar \ff$ as $r \lar 1$, so that \rrr{re2} holds. On the other hand, if $\sum_{n=1}^{\ff} |a_n|^2 < \ff$ then it is known that $f$ can be extended to a radial limit function $f(e^{i\th})$ on $\dd \DD$ for a.e. $\th$(see \cite[Thm. 17.10]{rud} or \cite[Sec. 6.5]{durBM}) such that $\lim_{r \nearrow 1} \int_{0}^{2\pi} |f(re^{i\th})-f(e^{i\th})|^2d\th = 0$ and $\sum_{n=1}^{\ff} |a_n|^2 = \frac{1}{2\pi}\int_{0}^{2\pi} |f(e^{i\th})|^2 d\th$. The result follows. An alternate proof of the lemma can be given using Green's function $G_r(z,0) = \frac{1}{\pi} \log\frac{r}{|z|}$, which satisfies the following "expected occupation time" formula(see \cite[Sec.'s 1.8,1.9]{durBM})

\begin{equation} \label{}
E\int_{0}^{\tau(r\DD,0)} g(B_t)dt = \int_{r\DD}g(z) G_r(z,0) dA(z)
\end{equation}

for any measurable, nonnegative function $g$. We have for $r \leq 1$

\begin{equation} \label{}
\begin{split}
E \int_{0}^{\tau(r\DD,0)} |f'(B_s)|^2 ds & = \int_{r\DD} |f'(z)|^2 G_r(z,0)dA(z) \\
& = \frac{1}{\pi}\int_{0}^{r} s \log \frac{r}{s} \int_{0}^{2\pi} |f'(se^{i\th})|^2 d\th ds \\
& = 2 \int_{0}^{r} s \log \frac{r}{s} \sum_{n=1}^{\ff} n^2 |a_n|^2 s^{2n-2} ds \\
& = 2 \sum_{n=1}^{\ff} n^2 |a_n|^2 \int_{0}^{r} s^{2n-1} \log\frac{r}{s} ds \\
& = 2 \sum_{n=1}^{\ff} n^2 |a_n|^2 r^{2n} \int_{0}^{1} s^{2n-1} \log \frac{1}{s} ds \\
& = 2 \sum_{n=1}^{\ff} n^2 |a_n|^2 r^{2n} (\frac{1}{4n^2}) = \frac{1}{2}\sum_{n=1}^{\ff} |a_n|^2 r^{2n}
\end{split}
\end{equation} \qed

Both theorems now follow. Theorem \ref{kris} and the upper bound in Theorem \ref{greg} are immediate from de Branges' Theorem because the coefficients of the Taylor series of any $f \in \SS$ are dominated by those of $K(z)$. The lower bound in Theorem \ref{greg} is also clear since $I(z)$ is the only function in $\SS$ with $a_n = 0$ for $n \geq 2$.

\vski

{\bf Remarks:} The expression \rrr{re2} is half of the square of the  $H^2$ Hardy norm of $f$, written $||f||^2_{H^2}$; see \cite{rud} for more details. In \cite{burk}, the $H^p$ norm of $f$ is used to study the finiteness of the $p/2$ moment of $\tau(f(\DD),f(0))$. The upper bound of Theorem \ref{greg} can also be taken to be the special case $p=2$ of the following result of Baernstein in \cite{burn}.

\begin{theorem} \label{just}
If $f\in \SS, f \notin \KK$, $r \in (0,1)$ then

\begin{equation} \label{}
\int_{0}^{2\pi} |f(re^{i\th})|^p d\th < \int_{0}^{2\pi} |K(re^{i\th})|^p d\th
\end{equation}

for any $p \in (0,\ff)$.
\end{theorem}

The proof of Theorem \ref{just} given in \cite{burn}, however, is quite involved, so it seemed as well to give a simple proof based on de Branges' Theorem(which was not available when \cite{burn} was published), as we have done. We remark that the fact that Baernstein's result can be related to the exit time of Brownian motion was first noted by Betsakos in \cite{bets2}. In light of the connection given by Burkholder in \cite{burk} between Hardy norms of analytic functions and moments of the exit times of Brownian motion, it is tempting to conjecture that the following is a consequence of Theorem \ref{just}.

\begin{conjecture} \label{}
If $f\in \SS, f \notin \KK$, $r \in (0,1)$, then
\begin{equation} \label{}
E_{0}[\tau(f(r\DD))^p] < E_{0}[\tau(K(r\DD))^p]
\end{equation}

for any $p \in (0,\ff)$. The same holds when $r=1$ and $p \in (0,1/4)$.
\end{conjecture}

\vski

Unfortunately, the inequalities in \cite{burk} do not seem to be sufficient to prove this result.

\section{Further consequences of Lemma \ref{bigguy}} \label{cons}

Lemma \ref{bigguy}, if noticed before, does not seem to have been used in calculating the exit times of domains. In this section we give a few examples where it can be applied with relative ease. We also point out how the lemma can be used in a different application, namely, the summation of certain series. That is, Lemma \ref{bigguy} expresses the expected exit time as a sum. If we have a different way of evaluating the expected exit time, then we have obtained a value for the sum. The following derivation of a result due to Euler is perhaps the simplest example of this idea; see Example \ref{square} below for a more complex instance.

\begin{prop} \label{}
\begin{equation} \label{all}
\sum_{n=1}^\ff \frac{1}{n^2} = \frac{\pi^2}{6}
\end{equation}
\end{prop}

{\bf Proof:} We will prove the corresponding statement obtained by summing only the odd terms, namely

\begin{equation} \label{odd}
\sum_{n=1}^\ff \frac{1}{(2n-1)^2} = \frac{\pi^2}{8}
\end{equation}

This is equivalent to \rrr{all}, as

\begin{equation} \label{}
\sum_{n=1}^\ff \frac{1}{(2n-1)^2} = \sum_{n=1}^\ff \frac{1}{n^2} - \sum_{n=1}^\ff \frac{1}{(2n)^2} = (1-\frac{1}{4})\sum_{n=1}^\ff \frac{1}{n^2} = \frac{3}{4}\sum_{n=1}^\ff \frac{1}{n^2}
\end{equation}

Let $W_t$ be a one-dimensional Brownian motion with $W_0=0$ a.s., and let $T= \inf_{t>0} \{|W_t| = \frac{\pi}{4}\}$. We will calculate $E[T]$ in two different ways. The first way is quite standard(see, for example, \cite[Ex. 7.5]{fima}). We apply the Optional Stopping Theorem to the martingale $W_t^2 - t$ to obtain

\begin{equation} \label{true}
E[T] = E[W_T^2] = \frac{\pi^2}{16}
\end{equation}

We now calculate $E[T]$ using Lemma \ref{bigguy}. $W_t$ may be taken to be the real part of our two dimensional Brownian motion $B_t$, and it is therefore clear that $E[T] = E_0[\tau(U)]$, where $U = \{ \frac{-\pi}{4} < \mbox{Re } z < \frac{\pi}{4} \}$. We need then to find a conformal map $f(z)$ mapping $\DD$ onto $U$ with $f(0) = 0$. The function

\be
\tan z = \frac{\sin z}{\cos z} = -i \frac{e^{2iz} - 1}{e^{2iz} + 1}
\ee

maps $U$ conformally to $\DD$. This can be seen by noting that the function $x+iy \lar e^{2i(x+iy)} = e^{-2y + 2ix}$ maps $U$ conformally to $\{ \mbox{Re } z > 0\}$, and then that the M\"obius transformation $z \lar -i(\frac{z-1}{z+1})$ maps $\{ \mbox{Re } z > 0\}$ conformally to $\DD$. We conclude that the principal branch of $\tan^{-1}z$ maps $\DD$ conformally to $U$. $\tan^{-1}z$ admits the Taylor series expansion

\begin{equation} \label{}
\tan^{-1}z = z - \frac{z^3}{3} + \frac{z^5}{5} - \ldots = \sum_{n=1}^{\ff} \frac{(-1)^{n+1} z^{2n-1}}{2n-1}
\end{equation}

Thus, by Lemma \ref{bigguy},

\begin{equation} \label{doubletrue}
E[T] = \frac{1}{2} \sum_{n=1}^{\ff} \frac{1}{(2n-1)^2}
\end{equation}

Equating \rrr{true} and \rrr{doubletrue} yields \rrr{odd}. \qed

Here are some further examples.

\vski

\ccases{disc}: Consider the disc $r\DD$ for any $r > 0$, and let $a \in r\DD$. It may be checked that

\begin{equation} \label{}
f(z) = r \frac{z+\frac{a}{r}}{1 + \frac{\bar{a}}{r} z} = r \frac{r z+a}{r + \bar{a} z}
\end{equation}
is a conformal map sending $\DD$ to $r\DD$ and $0$ to $a$. To find the power series for $f$ we expand

\begin{equation} \label{}
\begin{split}
r \frac{z+\frac{a}{r}}{1 + \frac{\bar{a}}{r} z} & = r (z+\frac{a}{r})(1-\frac{\bar{a}z}{r}+ (\frac{\bar{a}z}{r})^2 - \ldots) \\
& = a + (r^2 - |a|^2) \sum_{n=1}^{\ff} \frac{(-1)^{n-1} \bar{a}^{n-1}z^n}{r^{n}}
\end{split}
\end{equation}

Lemma \ref{bigguy} now gives

\begin{equation*} \label{}
E_a[\tau(r\DD)] = \frac{1}{2}(r^2 - |a|^2)^2 \sum_{n=1}^{\ff} \frac{|a|^{2n-2}}{r^{2n}} = \frac{(r^2 - |a|^2)^2}{2r^2} \frac{1}{1-\frac{|a|^2}{r^2}} = \frac{1}{2} (r^2 - |a|^2)
\end{equation*}

This is well known; see for example \cite[Ex. 7.4.2]{ox}.
\vski

\ccases{hp}: The function $f(z)=\frac{z}{1-z} = z + z^2 + z^3 + \ldots$ maps $\DD$ conformally to $\{\mbox{Re }z> -\frac{1}{2}\}$. Lemma \ref{bigguy} gives $E_0[\tau(\{\mbox{Re }z> -\frac{1}{2}\})] = \ff$, and it is easy to see that $E_z[\tau(U)] = \ff$ whenever $U$ is a half-plane and $z \in U$. Letting $W_t = \mbox{Re }B_t$ we recover a known result on one-dimensional Brownian motion, that if $T=\inf_{t>0}\{W_t = a\}$ with $a \neq 0$ then $E[T] = \ff$.

\vski

%\ccases{skip}: The function $f(z) = \frac{iz}{1-z^2} = i(z + z^3 + z^5 + \ldots)$ maps $\DD$ conformally onto $U=\CC \backslash (-\ff,1/2]\cup [1/2,\ff)$, and we see $E_z[\tau(U)] = \ff$ for $z \in U$. This also can be deduced easily from Example \ref{hp}.

%\vski

\ccases{cardiod}: Let $U$ be the cardioid with boundary defined by the polar equation $r = 2(1+\cos \th)$.

\hspace{1.2in} \includegraphics[width=70mm,height=70mm]{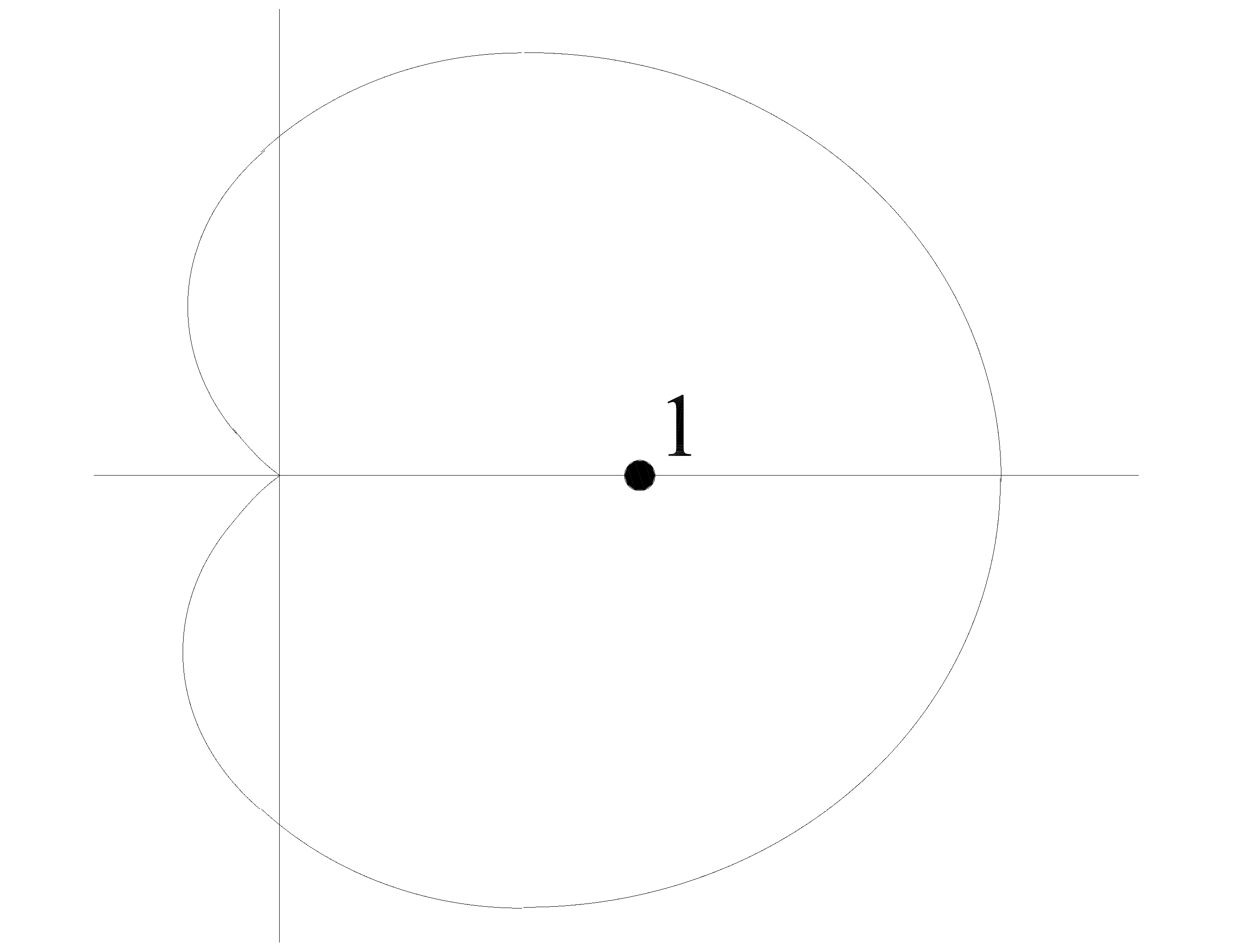}

This is the conformal image under $z^2$ of the disc $\{|z-1|=1\}$. The conformal map from $\DD$ to $U$ mapping $0$ to $1$ is given by $f(z)=(z+1)^2 = 1+2z +z^2$. Applying Lemma \ref{bigguy} we get $E_1[\tau(U)] = \frac{5}{2}$.

\vski

\ccases{logimage}: Let $\ga$ denote the curve in $\RR^2$ defined by $e^x = 2 \cos y$ for $\frac{-\pi}{2} < y < \frac{\pi}{2}$.

\hspace{1.2in} \includegraphics[width=100mm,height=70mm]{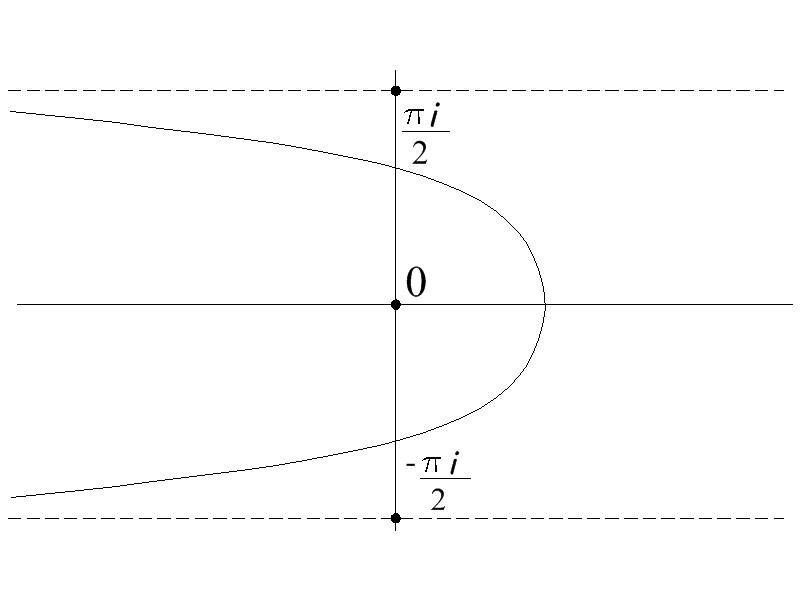}

This curve has been referred to as the "catenary of equal resistance". Let $U$ be the component of $\ga^c$ which contains 0. Then $U$ is the conformal image of $\DD$ under the map

\be
f(z) = \log(1+z) = z -\frac{z^2}{2} + \frac{z^3}{3} - \frac{z^4}{4} + \ldots
\ee

Applying Lemma \ref{bigguy} we obtain $E_0[\tau(U)] = \frac{1}{2}\sum_{n=1}^{\ff} \frac{1}{n^2} = \frac{\pi^2}{12}$.

\vski

\ccases{wedge}: Let $\HH = \{ \mbox{Re }z > 0\}$ be the right half-plane\footnotemark
\footnotetext{This is at odds with standard practice in complex analysis, where usually $\HH = \{ \mbox{Im }z > 0\}$ is the upper half-plane. It is convenient for our purposes, however.}
, and let $\HH^p = \{|Arg(z)| < \frac{\pi p}{2}\}$ for $p \leq 1$, where $Arg(z)$ is the principal branch of the argument function taking values in $(-\pi,\pi]$. $\HH^p$ is the infinite wedge of width $\pi p$ centered at the positive real axis, and $\HH^1 = \HH$. It was shown in $\cite{spitz}$, among other things, that $E_1[\tau(\HH^p)]< \ff$ if $p < \frac{1}{2}$, and $E_1[\tau(\HH^p)]= \ff$ if $p \geq \frac{1}{2}$. We will derive this result and find bounds for $E_1[\tau(\HH^p)]$ using Lemma \ref{bigguy}. The conformal map from $\DD$ to $\HH^p$ is $g(z) = \frac{(1+z)^p}{(1-z)^p}$, but the Taylor series expansion for $g$ appears to be unwieldy for arbitrary $p$, so we will simplify somewhat. Define the principal branch of $z^p$ on $\HH$, where $Arg(z^p) = p Arg(z)$. Then the inverse $z^{1/p}$ is well defined on $\HH^p$. Let $\tilde{\HH}^p = \{Re(z^{1/p})>1/2\}$. $\HH^p$ is then the image of $\{\mbox{Re }z > \frac{1}{2}\}$ under $z^p$. The relationship between $\HH^p$ and $\tilde{\HH}^p$ is shown below.

\hspace{1.2in} \includegraphics[width=70mm,height=50mm]{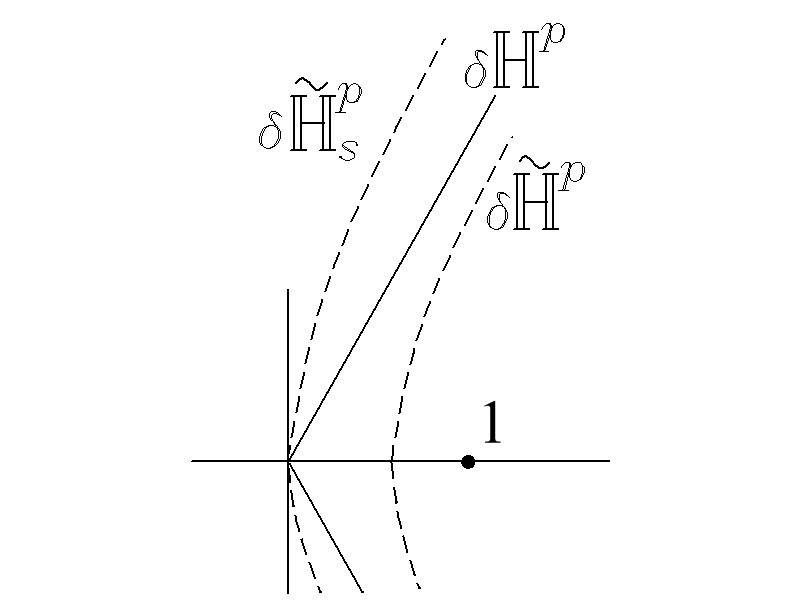}

It is clear that $\tilde{\HH}^p \subseteq \HH^p$, and it is also easy to check that $\HH^p \subseteq \tilde{\HH}^p_s := \{\frac{2^p-1}{2^p} z+\frac{1}{2^p} \in \tilde{\HH}^p\}$. The conformal map from $\DD$ to $\tilde{\HH}^p$ is given by $f(z) = \frac{1}{(1-z)^p}$, and the conformal map from $\DD$ to $\tilde{\HH}^p_s$ is given by $\frac{2^p}{2^p-1}(f(z) - 1/2^p)$. We have

\begin{equation} \label{}
\frac{d^n}{dz^n} \frac{1}{(1-z)^p} \Big| _{z=0} = p(p+1) \ldots (p+n-1)
\end{equation}

The expression $(p)_n := p(p+1) \ldots (p+n-1)$ with $(p)_0 = 1$ is known as the {\it Pochhammer symbol}, and we see that $f(z) = \sum_{n=0}^\ff \frac{(p)_n z^n}{n!}$. Applying Lemma \ref{bigguy} we obtain

\begin{equation} \label{}
E_1[\tau(\tilde\HH^p)] = \frac{1}{2}\sum_{n=1}^{\ff} \frac{(p)_n^2}{(n!)^2}
\end{equation}

This series can be evaluated explicitly. Using the notation in \cite{luke} it is given by $_2 F_1(p,p;1;1)-1$, where $_2 F_1$ refers to the hypergeometric function. Applying Euler's integral formula \cite[Sec. 3.6]{luke} and using the standard definitions for the $\Be$ and $\Ga$ functions, we calculate

\begin{equation} \label{}
E_1[\tau(\tilde\HH^p)] = \frac{1}{2}(_2F_1(p,p;1;1)-1) = \frac{1}{2}\Big(\frac{\Ga(1)\Be(p,1-2p)}{\Ga(1-p)\Ga(p)} - 1 \Big)
\end{equation}

Similarly, we have

\begin{equation} \label{}
E_1[\tau(\tilde\HH^p_s)] = \frac{2^{2p}}{2(2^p-1)^2} \Big(\frac{\Ga(1)\Be(p,1-2p)}{\Ga(1-p)\Ga(p)} - 1 \Big)
\end{equation}

We obtain

\begin{equation} \label{}
\frac{2^{2p}}{2(2^p-1)^2} \Big(\frac{\Ga(1)\Be(p,1-2p)}{\Ga(1-p)\Ga(p)} - 1 \Big) < E_1[\tau(\HH^p)] < \frac{1}{2}\Big(\frac{\Ga(1)\Be(p,1-2p)}{\Ga(1-p)\Ga(p)} - 1 \Big)
\end{equation}

This is finite for $p<\frac{1}{2}$, but infinite for $p \geq \frac{1}{2}$, as the integral defining $\Be(p,1-2p)$ diverges at $p=\frac{1}{2}$. This proves the given statement. Note that \cite{burk} also used the Hardy norm of $f$ in a similar manner to deduce a more general result.

\vski

\ccases{square} Let $m$ be an integer greater than 2 and set $\om = e^{\frac{2 \pi i}{m}}$. Let $U_m$ be the regular $m$-gon with vertices at $1, \om, \om^2, \ldots , \om^{m-1}$. We will calculate $E_0[\tau(U_m)]$. Consider the Schwarz-Christoffel mapping given by

\begin{equation} \label{}
g(z) = \int_0^z \frac{d\ze}{(1-\ze)^{2/m}(1-\om \ze)^{2/m} \ldots (1-\om^{m-1} \ze)^{2/m}} = \int_0^z \frac{d\ze}{(1-\ze^m)^{2/m}}
\end{equation}

$g(z)$ is a conformal mapping of the unit disc onto an $m$-gon with all angles $\pi^{(m-2)/m}$(see \cite{tref} for details), and symmetry arguments show that the image is in fact a regular $m$-gon with vertices $W,\om W,\ldots, \om^{m-1}W$ for some $W>0$. The points $1,\om,\ldots,\om^{m-1}$ map to these vertices, and we can calculate

\begin{equation} \label{}
\begin{split}
W & = g(1) = \int_0^1 \frac{dx}{(1-x^m)^{2/m}} = \frac{1}{m} \int_0^1 u^{-(m-1)/m}(1-u)^{-2/m}du \\
& = \frac{1}{m} \beta(1/m,(m-2)/m)
\end{split}
\end{equation}

Setting $f(z) = \frac{1}{W} g(z)$, we obtain our conformal map from $\DD$ to $U_m$. Recall from Example \ref{wedge} that $\frac{1}{(1-z)^{2/m}} = \sum_{n=0}^\ff \frac{(2/m)_n z^n}{n!}$. We see that

\begin{equation} \label{}
\begin{split}
f(z)  & = \frac{1}{W} \int \sum_{n=0}^\ff \frac{(2/m)_n z^{mn}}{n!} = \frac{1}{W} \sum_{n=0}^\ff \frac{(2/m)_n z^{mn+1}}{n!(mn+1)} \\
& = \frac{z}{W} \sum_{n=0}^\ff \frac{(1/m)_n(2/m)_n z^{mn}}{((m+1)/m)_nn!} = \frac{z}{W} {}_2F_1(1/m,2/m;(m+1)/m;z^m)
\end{split}
\end{equation}

where we have used the identity $(1/m)_n (mn+1) = ((m+1)/m)_n$. Theorem \ref{bigguy} gives

\begin{equation} \label{rte}
\begin{split}
E_0[\tau(U_m)] & = \frac{1}{2W^2}\sum_{n=0}^\ff \frac{(2/m)^2_n (1/m)^2_n }{(n!)^2((m+1)/m)_n^2} \\
& = \frac{1}{2W^2} {}_4 F_3(1/m,1/m,2/m,2/m;(m+1)/m,(m+1)/m,1;1) \\
& = {}_4 F_3(1/m,1/m,2/m,2/m;(m+1)/m,(m+1)/m,1;1) \\
& \hspace{1.3in} \times \frac{m^2}{2\beta(1/m,(m-2)/m)^2}
\end{split}
\end{equation}

In the cases $m=3,4$ we may compare this with known results. For the equilateral triangle, applying \cite[Thm. 1]{ala} gives $E_0[\tau(U_3)] = 1/6$. Computer approximation shows agreement with \rrr{rte} evaluated at $m=3$. We arrive at the following value for the hypergeometric sum.

\begin{equation} \label{}
{}_4 F_3(1/3,1/3,2/3,2/3;4/3,4/3,1;1) = \frac{\beta(1/3,1/3)^2}{27}
\end{equation}

For the case $m=4$, the square, \rrr{rte} gives

\begin{equation} \label{run}
E_0[\tau(U)] = {}_4F_3(1/4,1/4,1/2,1/2;5/4,5/4,1;1) \times \frac{8}{\beta(1/4,1/2)^2} \approx .294685
\end{equation}

This agrees with the approximation given in \cite[Tab. 10]{helm}\footnotemark \footnotetext{The value in \cite[Tab. 10]{helm} must be doubled, as the calculations there were for a square with unit side length rather than our normalization.}. This approximation was based on an explicit expression obtained from \cite{knight}\footnotemark \footnotetext{There is a misprint in \cite{knight}; the correct formula is given in \cite{helm}.}, and equating that expression with \rrr{run} we obtain the following strange identity.

\begin{equation} \label{}
\begin{split}
{}_4F_3(1/4,& 1/4,1/2,1/2;5/4,5/4,1;1) \times \frac{1}{\beta(1/4,1/2)^2} \\ & = \frac{8}{\pi^4} \sum_{n=1}^\ff \sum_{m=1}^{\ff} \frac{(-1)^{m+n}}{(2m-1)(2n-1)((2m-1)^2+(2n-1)^2)}
\end{split}
\end{equation}

The result \rrr{rte} for $m \geq 5$ may be new.

\section{Extension to arbitrary domains}

Although not the main focus of the paper, we would be remiss if we failed to observe that Lemma \ref{bigguy} can be extended to some analytic functions which are not conformal. Examining the proof of the theorem should reveal that the hypothesis of injectivity is not necessary for the statement to hold. Instead, the important property of conformal maps which was used is that $f(B_t)$ leaves $f(\DD)$ at time $\tau(\DD, 0)$. This is not true of arbitrary analytic maps. For example, let $f$ map $\DD$ conformally to the rectangle $V=\{-1 < \mbox{Re }z < 1, -2\pi < \mbox{Im }z < 2\pi\}$. Then $g=e^f$ maps $\DD$ onto the annulus $\{\frac{1}{e} < |z| < e\}$. The boundary segments $\ell_1=\{-1 < \mbox{Re }z < 1,\mbox{Im }z = 2\pi\}$ and $\ell_2=\{-1 < \mbox{Re }z < 1,\mbox{Im }z = -2\pi\}$ are mapped by $e^z$ to the interior segment $\{\frac{1}{e} < \mbox{Re }z < e,\mbox{Im }z = 0\}$. We see that for any Brownian path $B_t(\om)$ such that $B_{\tau(\DD,z)}(\om) \in f^{-1}(\ell_1 \cup \ell_2)$ we have $g(B_{\tau(\DD,z)}(\om)) \in \{\frac{1}{e} < \mbox{Re }z < e,\mbox{Im }z = 0\}$. Since $g(B_t)$ does not leave $g(\DD)$ with probability 1 at time $\tau(\DD, 0)$, Theorem \ref{bigguy} will fail to hold.

\vski

With this in mind, let us define an analytic function $f$ on $\DD$ to be $B-proper$ if a.s. $f(B_t)$ leaves every compact subset of $f(\DD)$ as $t$ increases to $\tau(\DD,0)$. Let a domain be called $B-proper$ if it is the image of a B-proper map on $\DD$. The reason for this terminology is that analytic functions with the property that $f(z_n)$ leaves every compact set for every sequence $\{z_n\}$ which approaches $\dd \DD$ are commonly referred to as $proper$. It is easy to see that every conformal map is proper. It is also clear that every proper function is B-proper, but the converse is not true, as the example given below shows. With the same proof as Lemma \ref{bigguy} we have the following.

\begin{lemma} \label{bigguy2}
Suppose $f(z) = \sum_{n=0}^{\ff} a_n z^n$ is B-proper on $\DD$. Then
\begin{equation} \label{re3}
E_{f(0)}[\tau(f(\DD))] = \frac{1}{2}\sum_{n=1}^{\ff} |a_n|^2
\end{equation}
\end{lemma}

We now give an example showing that a function can fail to be proper but still be B-proper. Let $f(z) = e^{tan^{-1}z}$. $f$ maps $\DD$ conformally to the annulus $\{e^{\frac{-\pi}{4}} < |z| < e^{\frac{\pi}{4}}\}$. If $z_n$ is a sequence in $\DD$ approaching $i$ or $-i$ along the imaginary axis, $f(z_n)$ simply cycles around $A$ on the circle $\{ |z|=1 \}$. $f$ is, however, B-proper, since $f$ extends continuously to map $\{z=e^{i\th}; -\frac{\pi}{2}<\th <\frac{\pi}{2}\}$ to $\{ |z| = e^{\frac{\pi}{4}}\}$ and $\{z=e^{i\th}; \frac{\pi}{2}<\th <\frac{3\pi}{2}\}$ to $\{ |z| = e^{\frac{-\pi}{4}}\}$. We see that $f(B_t)$ leaves $f(\DD)$ with probability 1 as $t$ increases to $\tau(\DD)$.

%Note that this implies the following result, which is trivial for conformal maps but less so for B-proper or even proper ones.

%\begin{corollary} \label{}
%Suppose $f$ and $g$ are B-proper maps on $\DD$, with $f(\DD)=g(\DD)$ and $f(0)=g(0)$. Then $||f||_{H^2} = ||g||_{H^2}$.
%\end{corollary}

\section{Acknowledgements}

I would like to thank George Markowsky and Kais Hamza for helpful conversations. I am also grateful for support from Australian Research Council Grant DP0988483.

%\bibliography{CFR}

\def\noopsort#1{} \def\printfirst#1#2{#1} \def\singleletter#1{#1}
   \def\switchargs#1#2{#2#1} \def\bibsameauth{\leavevmode\vrule height .1ex
   depth 0pt width 2.3em\relax\,}
\makeatletter \renewcommand{\@biblabel}[1]{\hfill#1.}\makeatother

\bibliographystyle{alpha}
\bibliography{CABMbib}

\end{document}